\documentclass[11pt]{article}
\usepackage{latexsym,amscd,amsmath,amssymb}
\title{\bf Biholomorphic Convex Mappings of Order $\alpha$ on the Unit Ball in Hilbert
Spaces}
\author{Ming-Sheng Liu\thanks{Corresponding author.} \\
School of Mathematical Sciences, South China Normal University,\\
Guangzhou, 510631 Guangdong, P.R.China\\
Email:liumsh@scnu.edu.cn\\
Yu-Can Zhu \\
Department of Mathematics, Fuzhou University, Fuzhou,350002 Fujian,
P.R.China \\{E-mail:zhuyucan@fzu.edu.cn}}
\date{}
\begin{document}
\maketitle \large
\def\no{\noindent}
\def\dl{\displaystyle\lim}
\def\ds{\displaystyle\sum}
\def\di{\displaystyle\int}
\newcommand{\f}{\displaystyle\frac}
\newcommand{\il}{\int\limits}
\newcommand{\s}{\sum\limits}
\newcommand{\ba}{\begin{array}}
\newcommand{\ea}{\end{array}}
\newcommand{\ld}{\left\{}
\newcommand{\rd}{\right\}}
\newcommand{\lz}{\left[}
\newcommand{\rz}{\right]}
\newcommand{\lx}{\left(}
\newcommand{\rx}{\right)}
\newcommand{\og}{\omega}
\newcommand{\ra}{\rightarrow}
\newcommand{\he}[2]{\sum\limits_{#1}^{#2}}
\newcommand{\fd}[4]{\left \{ \begin{array}{ll}
#1&#2\\#3&#4\end{array}\right.}
\renewcommand{\thefootnote}{\fnsymbol{footnote}}

\begin{center}
\begin{minipage}{14cm}

{\bf Abstract.} In this paper, we first introduce the concept of biholomorphic convex
mapping of order $\alpha$ on the unit
ball $B$ in a complex Hilbert space $X$. Next we provide some sufficient
conditions that a locally biholomorphic mapping $f$ is a biholomorphic convex
mapping of order $\alpha$ and give an Alexander's theorem between the subclass of convex mappings and the subclass of starlike mappings on $B$ in Hilbert space. We also obtain the order of starlikeness of biholomorphic convex mappings
 of order $\alpha$ on $B$ in Hilbert spaces. Finally, we construct some concrete examples of biholomorphic convex mappings of order $\alpha$ on $B$ in Hilbert spaces by means of a linear operator.

{\bf 2000 MR Subject Classification} \ 32H02, 30C45
\end{minipage}
\end{center}

\section{ Introduction}
\mbox{}\indent The holomorphic functions of one complex variable
which map the unit disk onto starlike or convex domains have been
extensively studied. These functions are easily characterized by
simple analytic or geometric conditions. In moving to higher
dimensions, several difficulties arise. In the case of one complex
variable, the following well known theorems had been established(cf. \cite{gk}).
\par
{\bf Theorem A}\, Suppose that $\alpha\in [0, 1)$ and
$f(z)=z+\sum\limits_{n=2}^\infty a_nz^n$ is a holomorphic function
on the unit disk $U=\{z\in\mathbb{C}: |z|<1\}$ in the complex plane $\mathbb{C}$.

(1)\, If $\sum\limits_{n=2}^\infty n^2 |a_n|\leq 1$, then $f$ is a convex
function in the unit disk $U$.

(2)\, If $\sum\limits_{n=2}^\infty n(n-\alpha )|a_n|\leq 1-\alpha$, then $f$ is a convex
function of order $\alpha$ in $U$.

{\bf Theorem B}\, Suppose that $\alpha\in [0, 1)$ and
$f(z)=z+\sum\limits_{n=2}^\infty a_nz^n$ is a holomorphic function
on the unit disk $U$ in the complex plane $\mathbb{C}$.

(1)\, If $\sum\limits_{n=2}^\infty n|a_n|\leq 1$, then $f$ is a starlike
function in the unit disk $U$.

(2)\, If $\sum\limits_{n=2}^\infty (n-\alpha )|a_n|\leq 1-\alpha$, then $f$ is a starlike
function of order $\alpha$ in $U$.

Roper and Suffridge established the n-dimensional version of
Theorem A(1), and we\cite{lz3} established the n-dimensional version of
Theorem B(1)(2) as follows.

{\bf Theorem C}(Roper and Suffridge\cite{rs})\quad Let $f(z)=z + \sum\limits_{k=2}^\infty
A_k(z^k)$ be a holomorphic mapping on the unit ball $B_2^n$. If
$\sum\limits_{k=2}^\infty k^2\|A_k\|\leq 1$,
then $f(z)$ is a convex mapping on $B_2^n$.

{\bf Theorem D}(Liu and Zhu\cite{lz3}) Suppose that $\alpha\in [0, 1)$. Let $f(z)=z +
\sum\limits_{k=2}^\infty A_k(z^k)$ be a holomorphic mapping on the unit ball
$B$ in Hilbert space. If $
\sum\limits_{k=2}^\infty (k-\alpha)\|A_k\|\leq A(\alpha)$,
where $A(\alpha)$ is defined by
$$
A(\alpha)= \left\{
\begin{array}{lll}
\frac{(2-\alpha)\sqrt{1-2\alpha}} {\sqrt{5-2\alpha}},&
0\leq\alpha\leq \frac{1}{4},\\
\frac{(2-\alpha)(1-\alpha)}{2+\alpha},&
\frac{1}{4}<\alpha\leq\frac{2}{5},\\
\alpha, &\frac{2}{5}<\alpha<\frac{1}{2},\\
1-\alpha, &\frac{1}{2}\leq\alpha <1.
\end{array}
\right.$$
Then $f(z)$ is a starlike mapping of order $\alpha$ on $B$  in Hilbert space.

A problem is naturally proposed: can we establish the n-dimensional
version for Theorem A(2)?

\section{Preliminaries}
\mbox{}\indent
In order to state and prove our main results, we recall some definitions
and notations. Suppose that $X$ is a complex  Hilbert space with
inner product $\langle \cdot, \cdot\rangle$ and norm
$\|\cdot\|=\sqrt{\langle \cdot, \cdot\rangle}$, and $G$ is a
domain in $X$. A mapping $f:G\to X$ is said to be holomorphic on
$G$, if for any $z\in G$, there exists a linear operator
$Df(z):X\to X$ such that
$$
\lim_{h\to 0}\frac{\| f(z+h)-f(z)-Df(z)h\|}{\| h \|}=0.
$$
The linear operator $Df(z)$ is called the Fr$\acute{e}$chet
derivative of $f$ at $z\in G$.

If $f$ is holomorphic on $G$, then for every $k=1, 2, \cdots $,
and every $z_0\in G$,  there is a bounded symmetric $k-$ linear
operator $D^{k}f(z_0):X\times X\times \cdots \times X\to X$ such
that
$$
f(z)=\he {k=0}{+\infty }\frac {1}{k!}D^{k}f(z_0)((z-z_0)^{k})
$$
for all $z$ in some neighborhood of $z_0$, where
$D^{0}f(z_0)((z-z_0)^{0})=f(z_0)$ and
$$
D^{k}f(z_0)((z-z_0)^{k})=D^{k}f(z_0)(z-z_0,z-z_0,\cdots , z-z_0)
$$
for $k\geq 1$.

A mapping $f:G\to X$ is said to be biholomorphic on $G$ if $f$ is
holomorphic on $G$, $f(G)$ is a domain, and the inverse $f^{-1}$
exists and is holomorphic on $f(G)$. A mapping $f:G\to X$ is said
to be locally biholomorphic on $G$, if for any $z\in G$, there
exists a neighborhood $U$ of $z$ such that $f|_U$ is biholomorphic
on $U$. Then $f$ is locally biholomorphic on $G$ if and only if
its Fr\'{e}chet derivative $Df(z)$ has a bounded inverse at each
$z\in G$.

The unit ball in $X$ is $B=\{z\in X: \|z\|<1\}$. Let $N(B)$ denote the class of all local biholomorphic
mappings $f:B\rightarrow X $ such that $f(0)=0,
Df(0)=I$, where $I$ is the identity operator in $X$. A biholomorphic mapping $f:B\to X$ is called a biholomorphic starlike mapping  if $tf(B)\subset f(B)$ for $0\leq t\leq1$ with $f(0)=0$. Let $S^{\ast}(B)$ be the subclass of $N(B)$
consisting of starlike mappings on $B$. Then $f\in S^{\ast}(B)$ if and only if $f$ is locally
biholomorphic such that
$$
\mbox{Re}\langle Df(z)^{-1}f(z),z\rangle>0
$$
for all $z\in B\backslash\{0\}$(cf.\cite{gs, gk, k}).

A mapping $f\in N(B)$ is called starlike of
order $\alpha\in (0,1)$ on $B$ if
$$
\bigg|\langle Df(z)^{-1}f(z),
z\rangle-\frac{1}{2\alpha }\|z\|^2\bigg|<\frac{1}{2\alpha }\|z\|^2 \quad \mbox{for
all }\quad z\in B\backslash\{0\},
$$
Let $S^{\ast}(B, \alpha)$ denote the class of starlike mappings of
order $\alpha$ on $B$ for $0<\alpha<1$ and let $S^{\ast}(B, 0)\equiv S^{\ast}(B)$. It is obvious that $S^{\ast}(B, \alpha)\subset S^{\ast}(B)$ for $0\leq\alpha<1$.

A biholomorphic mapping $f:B\to X$ is said
biholomorphic convex mapping if
$$
(1-t)f(z_1)+tf(z_2)\in f(B)
$$
for all $z_1, z_2\in B$ and $0\leq t\leq 1$. The class of all
biholomorphic convex mappings on $B$ with $f(0)=0, Df(0)=I$ is
denoted by $K(B)$.

We\cite{zl2} obtained a necessary and sufficient condition that a
locally biholomorphic mapping was a biholomorphic convex mapping on
B in the Hilbert space $X$ as follows.
\par
{\bf Theorem A(Zhu and Liu\cite{zl2}).} Let $f:B\to X$ be a locally biholomorphic mapping.
Then $f$ is a biholomorphic convex mapping on $B$ if and only if
\begin{eqnarray}
\|x\|^2-\mbox{Re}\Big\langle Df(z)^{-1}D^2f(z)(x, x),
z\Big\rangle>0\label{a01}
\end{eqnarray}
for $z\in B\setminus\{0\}$ and $x\in X\setminus\{0\}$ with
$\mbox{Re}\langle x, z\rangle=0$.

{\bf Remark 1.} Theorem A had improved a result
of Hamada and Kohr\cite{hk1}. Setting $X=C^n$ in Theorem A, we also obtain Theorem 2 in
\cite{gwy}.
\par
{\bf Corollary 1.} Let $0\leq \alpha<1$ and $f:B\to X$ be a locally biholomorphic
mapping. If $f$ satisfies the following inequality
\begin{eqnarray}
\|x\|^2-\mbox{Re}\Big\langle Df(z)^{-1}D^2f(z)(x, x),
z\Big\rangle>\alpha\cdot\|x\|^2\label{a02}
\end{eqnarray}
for $z\in B\setminus\{0\}$ and $x\in X\setminus\{0\}$ with
$\mbox{Re}\langle x, z\rangle=0$. Then $f$ is a biholomorphic
convex mapping on $B$.

We call such mapping $f$, which satisfies the hypothesis of Corollary 1, a biholomorphic
convex mapping of order $\alpha$ on $B$. We let $K(B, \alpha)$
denote the subclass of all biholomorphic convex mappings of order
$\alpha$ on $B$ with $f(0)=0, Df(0)=I$.

In this paper, we provide some sufficient conditions for biholomorphic convex mapping
of order $\alpha$ and an Alexander's theorem between the subclass of convex mappings and the subclass of starlike mappings on $B$ in Hilbert space. We also obtain the order of starlikeness of biholomorphic convex mappings
 of order $\alpha$ on $B$ in Hilbert spaces. Finally, we introduce a linear operator in purpose to
construct some concrete examples of biholomorphic convex mappings on
$B$ in Hilbert spaces. From these, we give some examples of
biholomorphic convex mappings on $B$ in Hilbert spaces.

\setcounter{equation}{0}
\section{Main Results}
\mbox{}\indent We first establish some sufficient conditions for biholomorphic convex
mapping of order $\alpha$ on $B$.

{\bf Theorem 1.} Suppose that $0\leq\alpha<1$ and $f:B\to X$ is a
locally biholomorphic mapping. If $f$ satisfies
\begin{eqnarray}
\Big\|Df(z)^{-1}D^2f(z)(x, x)\Big\|\leq 1-\alpha\label{a03}
\end{eqnarray}
for $z\in B\setminus\{0\}$ and $x\in X$ with $\|x\|=1$ and
$\mbox{Re}\langle x,z\rangle =0$, then $f$ is a biholomorphic convex
mapping of order $\alpha$ on $B$.

{\bf Proof.} Since $f:B\to X$ is a
locally biholomorphic mapping, for any
$z\in B\setminus\{0\}$ and $x\in X\setminus\{0\}$ with $\|x\|=1$ and
$\mbox{Re}\langle x,z\rangle =0$, we have
$$
\begin{array}{lll}
\|x\|^2-\mbox{Re}\Big\langle Df(z)^{-1}D^2f(z)(x, x),
z\Big\rangle&\geq&\|x\|^2-|\langle Df(z)^{-1}D^2f(z)(x, x),
z\rangle|\\
&\geq&\|x\|^2-\| Df(z)^{-1}D^2f(z)(x, x)\|\|z\|\\
&>&\|x\|^2-\| Df(z)^{-1}D^2f(z)(x, x)\|\\
&=&\|x\|^2-\| Df(z)^{-1}D^2f(z)(\frac{x}{\|x\|}, \frac{x}{\|x\|})\|\|x\|^2.
\end{array}
$$
\indent
Notice that $\|\frac{x}{\|x\|}\|=1$, we conclude from (\ref{a03}) that
$$
\begin{array}{lll}
\|x\|^2-\mbox{Re}\Big\langle Df(z)^{-1}D^2f(z)(x, x),
z\Big\rangle&>&\|x\|^2-\| Df(z)^{-1}D^2f(z)(\frac{x}{\|x\|}, \frac{x}{\|x\|})\|\|x\|^2\\
&\geq&\alpha\cdot\|x\|^2.
\end{array}
$$
Hence by Corollary 1, we obtain that
$f\in K(B, \alpha)$, and the proof is complete.

{\bf Corollary 2.} Suppose that $0\leq\alpha<1$ and $f:B\to X$ is a locally
biholomorphic mapping with $\|Df(z)-I\|\leq c<1$ for each $z\in
B$, where $I$ is the identity operator in $X$. If $f$ satisfies
$$
\|D^2f(z)(x, x)\|\leq (1-c)(1-\alpha)
$$
for all $x\in X$ with $\|x\|=1$ and $z\in B\setminus\{0\}$ such
that $\mbox{Re}\langle x, z\rangle =0$, then $f$ is a
biholomorphic convex mapping of order $\alpha$ on $B$.
\par
{\bf Proof.} Since $\|Df(z)-I\|\leq c<1$ for any $z\in B$, we
obtain that $Df(z)=I-(I-Df(z))$ is an invertible linear
operator(see \cite{tl}, P192), and
$$
\|Df(z)^{-1}\|\le\frac{1}{1-\|I-Df(z)\|} \le\frac{1}{1-c}
$$
for all $z\in B$.

Thus for any $x\in X$ with
$\|x\|=1$  and $z\in B\setminus\{0\}$ such that $\mbox{Re}\langle
x, z\rangle =0$, by the hypothesis of Corollary 2, we have
\begin{eqnarray*}
\|Df(z)^{-1}D^{2}f(z)(x, x)\|&\leq & \|Df(z)^{-1}\|\|D^2f(z)(x, x )\|\\
&\leq &\frac{1}{1-c}\cdot (1-c)(1-\alpha)=1-\alpha.
\end{eqnarray*}

Hence by Theorem 1, we obtain that $f$ is a biholomorphic convex
mapping of order $\alpha$ on $B$, and the proof is complete.
\par
{\bf Remark 2.} Setting $\alpha=0$ in Theorem 1, we get Corollary 1 in \cite{zl2}; Setting $\alpha=0$ in Corollary 2, we get Corollary 2 in \cite{zl2}.

{\bf Theorem 2.} Let $0\leq\alpha<1$ and $f(z)=z+\he{k=2}{+\infty}A_k(z^k):B\to X$ be a holomorphic mapping.
If $f$ satisfies
$\displaystyle\he{k=2}{+\infty}k(k-\alpha)\|A_k\|\leq 1-\alpha$,
then $f\in K(B, \alpha)$.

{\bf Proof.}\ By direct calculating the Fr$\acute{e}$chet
derivatives of $f(z)$, we obtain
\begin{eqnarray*}
&& Df(z)=I+\he{k=2}{+\infty}k A_k(z^{k-1}, \cdot),\\
&&D^2f(z)(x, x)=\he{k=2}{+\infty}k(k-1)A_k(z^{k-2},x^2)
\end{eqnarray*}
and
$$
\|Df(z)-I\|\leq\he{k=2}{+\infty}k\|A_k\|
\leq\frac{1}{2-\alpha}\he{k=2}{+\infty}k(k-\alpha)\|A_k\|\leq\frac{1-\alpha}{2-\alpha}<1
$$
for $z\in B, \, x\in X$.  Hence  we obtain that
$Df(z)=I-(I-Df(z))$ is an invertible linear operator(see
\cite{tl},P192), and
\begin{eqnarray*}
\|D^2f(z)(x,x)\|&\leq&\he{k=2}{+\infty}(k^2-k)\|A_k\|\|z\|^{k-2}\|\|x\|^{2}\\
&\leq& 1-\alpha +\alpha\he{k=2}{+\infty}k\|A_k\|-\he{k=2}{+\infty}k\|A_k\|\\
&=&(1-\he{k=2}{+\infty}k\|A_k\|)(1-\alpha) \label{a05}
\end{eqnarray*}
for $z\in B\setminus\{0\}$ and $ \|x\|=1$ with $\mbox{Re}\langle
x,z\rangle =0$. By Corollary 2 for
$c=\he{k=2}{+\infty}k\|A_k\|$, we obtain that
$f\in K(B, \alpha)$, and the proof is complete.

{\bf Remark 3.} Setting $X=\mathbb{C}^n, \alpha=0$ in Theorem 2, we may obtain
Theorem 2.1 in \cite{rs}. Our proof is more simple than theirs. Setting $X=\mathbb{C}$ in Theorem 2, we get
Theorem A(2).

{\bf Example 1.} Let $0\leq\alpha<1$ and $A$ be a symmetric bilinear operator from
$X\times X$ to $X$ with $\|A\|\leq\frac{1-\alpha}{4-2\alpha}$. If we let
$f(z)=z+A(z, z)$, then $f\in K(B, \alpha)$.

{\bf Proof.} Some straightforward computations yield the relations
$$
Df(z)=I+2A(z, \cdot),\quad D^2f(z)(x, y)=2A(x, y)
$$
for $z\in B, x, y\in X$. It implies
$$
Df(0)=I,\quad D^2f(0)(\cdot, \cdot)=2A(\cdot,
\cdot)\quad\mbox{and}\quad D^kf(0)=0
$$
for $k=3, 4, \cdots$. Hence we obtain
$$
\he{k=2}{+\infty}\frac{k(k-\alpha)\|D^kf(0)\|}{k!}=(2-\alpha)\|D^2f(0)\|=2(2-\alpha)\|A\|\leq 1-\alpha.
$$

By Theorem 2, we conclude that $f\in K(B, \alpha)$, and the proof is complete.

{\bf Example 2.} Let $0\leq\alpha<1$, $0<|a|\leq 1/2$ and $u\in X$ with $\|u\|=1$.
Then
$$
f(z)=z+a\langle z, u\rangle^2u\in K(B, \alpha)\Longleftrightarrow |a|\leq
\frac{1-\alpha}{4-2\alpha}.
$$

{\bf Proof.} Let $c=1-\frac{2|a|}{1-\alpha}$. If $0<|a|\leq \frac{1-\alpha}{4-2\alpha}$,  then we have
$\frac{1-\alpha}{2-\alpha}\leq c<1$ and $2|a|=(1-c)(1-\alpha)\leq c$. Short computations yield the
relations
\begin{eqnarray}
Df(z)=I+2a\langle z, u\rangle\langle\cdot, u\rangle u,\quad
D^2f(z)(x, x)=2a\langle x, u\rangle^2u.\label{a04}
\end{eqnarray}
It implies
$$
\|Df(z)-I\|\leq2|a|\|z\|<2|a|\leq c,\quad \|D^2f(z)(x,
x)\|\leq2|a|=(1-c)(1-\alpha)
$$
for all $x\in X$ with $\|x\|=1$ and $z\in B$ such that
$\mbox{Re}\langle x, z\rangle =0$.

By Corollary 2, we obtain that $f\in K(B, \alpha)$.

Conversely, we shall prove that $0<|a|\leq \frac{1-\alpha}{4-2\alpha}$ when $f\in K(B, \alpha)$.

If not, then $|a|>\frac{1-\alpha}{4-2\alpha}$. Let $\theta=\arg a, \, x=ie^{-i\theta}u$
and $z_0=-re^{-i\theta}u$ for $\frac{1-\alpha}{(4-2\alpha)|a|}<r<1$, where $u\in X$
with $\|u\|=1$. Then $\|x\|=1, z_0\in B\setminus\{0\}$ and
$\mbox{Re}\langle x, z_0\rangle=\mbox{Re}\{-ir\}=0$.

Some straightforward computations from (\ref{a04}) yield the
relations
\begin{eqnarray*}
&&Df(z_0)^{-1}=I+\frac{2|a|r}{1-2|a|r}\langle \cdot, u\rangle u,
\quad D^2f(z_0)(x, x)=-2|a|e^{-i\theta}u\\
&&Df(z_0)^{-1}D^2f(z_0)(x, x)=-\frac{2|a|e^{-i\theta}}{1-2|a|r}u.
\end{eqnarray*}
Hence we obtain
$$
\|x\|^2-\mbox{Re}\Big\langle Df(z_0)^{-1}D^2f(z_0)(x, x),
z_0\Big\rangle=\frac{1-4|a|r}{1-2|a|r}<\alpha.
$$
This contradicts (\ref{a02}), and the proof is complete.

Next, we provide an Alexander's theorem between the subclass of convex mappings and the subclass of starlike mappings on $B$ in Hilbert space. In the case of one complex variable, Alexander's theorem told us that $f(z)$ is a convex function on the unit disc
$ U  $ if and only if $zf^{\prime}(z)$ is a starlike
function on the unit disc $ U $. This theorem is no longer
true in several complex variables(see \cite{gk}). However, we have the following Alexander's theorem.

{\bf Theorem 3(Alexander's Theorem).}\,  Suppose that $0\leq\alpha<1$ and $A(\alpha)$ is defined by Theorem D. Let
$$
SK(B, \alpha)=\{f(z): f(z)=z + \sum\limits_{k=2}^\infty
A_k(z^k)\in H(B)\mbox{ such that } \he{k=2}{\infty}k(k-\alpha)\parallel A_k\parallel\leq A(\alpha)\},
$$
and
$$
SS^{*}(B, \alpha)=\{f(z): f(z)=z + \sum\limits_{k=2}^\infty
A_k(z^k)\in H(B)\mbox{ such that } \he{k=2}{\infty}(k-\alpha)\parallel A_k\parallel\leq A(\alpha)\}.
$$
Then $SK(B, \alpha)\subset K(B)$, $SS^{*}(B, \alpha)\subset S^{*}(B)$, and $f(z)\in SK(B, \alpha)$ if and only if
$Df(z)(z)\in SS^{*}(B, \alpha)$.

Notice that $A(\alpha)\leq 1-\alpha$ for $0\leq\alpha<1$, by applying Theorem D and Theorem 2, we can prove this theorem easily.

Now we establish a result on the order of starlikeness of function class $K(B, \alpha)$.

{\bf Theorem 4.} Suppose that $0\leq\alpha<1$. Then $
K(B, \alpha)\subset S^{*}(B, \beta)$, where
$$\beta=\beta(\alpha)=\frac{2\alpha-1+\sqrt{(2\alpha-1)^2+8}}{4}.$$

In order to prove the above theorem, we need the following lemmas.

{\bf Lemma 1.}(\cite{mm})\ Let $g(z)=a+a_1 z+\cdots$ is analytic in $U$ and $g(z)\not\equiv a$. If there exists $z_0\in U\backslash\{0\}$ such that
 $|g(z_0)|=\max\limits_{|z|\leq|z_0|}|g(z)|$, then there exists real number $t\geq \frac{|g(z_0)|-|a|}{|g(z_0)|+|a|}$ such that $z_0g^{'}(z_0)=t\ g(z_0)$.

{\bf Lemma 2.}\ Let $f\in N(B),\, 0<\rho<1$. If there exists $z_0\in B\backslash\{0\}$ such that
$$
{\rm Re}\ \frac{\|z\|^2}{\langle Df(z_0)^{-1}f(z_0), z_0\rangle}=\rho,
$$
and ${\rm Re}\ \frac{\|z\|^2}{\langle Df(z)^{-1}f(z), z\rangle}\geq\rho$ for all $\|z\|<\|z_0\|$.
Then there exist real numbers $\theta,\ t,\, m$ such that:

(i)\ $\langle h(z_0), z_0\rangle=\frac{\|z_0\|^2}{2\rho}\ (1+e^{i\theta}) $, where $h(z)=Df(z)^{-1}f(z)$ ;

(ii)\ $\langle Dh(z_0)(z_0), z_0\rangle=\frac{\|z_0\|^2}{2\rho}\ [1+(1+t)e^{i\theta}] $, where $t\geq \frac{1-|2\rho-1|}{1+|2\rho-1|}$ ;

(iii)\ $e^{i\theta}\overline{Dh(z_0)}(z_0)+e^{-i\theta}h(z_0)=m z_0$, where $m=\frac{2\cos\theta + 2 + t}{2\rho}$.

{\bf Proof.}\ Since $h:B\to X$ is a holomorphic mapping and $h(0)=0, Dh(0)=I$, and
\begin{equation}
\bigg|\langle h(z_0), z_0\rangle- \frac{\|z_0\|^2}{2\rho}\bigg|=\frac{\|z_0\|^2}{2\rho},
\label{a1}
\end{equation}
\begin{equation}
\bigg|\langle h(z), z\rangle- \frac{\|z\|^2}{2\rho}\bigg|\leq\frac{\|z\|^2}{2\rho},\quad \|z\|<\|z_0\|.
\label{a2}
\end{equation}

Let
$$\psi (\xi)=\frac{2\rho}{\|\xi z_0\|^2}\langle h(\xi z_0), \xi z_0\rangle - 1=\frac{2\rho}{\|z_0\|^2}\langle \frac{h(\xi z_0)}{\xi}, z_0\rangle - 1,
$$
then $\psi (\xi)$ is analytic in $\overline{U}$ and $|\psi (\xi)|\leq 1=|\psi (1)|$ for $\xi\in\overline{U}$, $\psi (0)=2\rho-1$. By Lemma 1, there is a real number $t\geq \frac{1-|2\rho-1|}{1+|2\rho-1|}$ such that $\psi^{'}(1)=t\psi(1)$.

Let $\psi(1)=e^{i\theta}$ for some real number $\theta$, then we obtain (i) holds, and
$$
\psi^{'}(1)=\frac{2\rho}{\|z_0\|^2}\langle Dh(z_0)(z_0)-h(z_0), z_0\rangle=te^{i\theta},
$$
which implies (ii) holds.

From (\ref{a1}) and (\ref{a2}), we obtain that
\begin{equation}
{\rm Re}[e^{-i\theta}\langle h(z), z\rangle ]\leq\frac{\|z_0\|^2}{2\rho}(1+\cos\theta)={\rm Re}[e^{-i\theta}\langle h(z_0), z_0\rangle ],\quad \|z\|<\|z_0\|.
\label{a3}
\end{equation}

Let $r=\|z_0\|, B(r)=\{z\in X: \|z\|<r\}$, then the tangent hyperplane of $\partial B(r)$ at $z_0$ is
$$
T_{z_0}(\partial B(r))=\{b\in X: \langle b, z_0\rangle=0\}.
$$

For any tangent vector $a\in T_{z_0}(\partial B(r))$ with $\|a\|=1$, set $\gamma (t)=\sqrt{1-t^2}\ z_0+t\ r a$, then $\|\gamma(t)\|=r$ for $t\in (-1, 1)$ and $\gamma(0)=z_0, \gamma'(0)=r a$. Let
$$
\varphi(t)={\rm Re}[e^{-i\theta}\langle h(\gamma (t)), \gamma (t)\rangle ].
$$

From (\ref{a3}), we obtain $\varphi(t)\leq \varphi(0)$ for $t\in (-1, 1)$, so that $\varphi(0)=\max\limits_{t\in (-1, 1)}\varphi(t)$. Hence
\begin{eqnarray*}
0=\varphi'(0)&=&{\rm Re}[e^{-i\theta}\langle D h(z_0)r a, z_0\rangle +  e^{-i\theta}\langle h(z_0), r a\rangle ]
= r {\rm Re}\langle v, a\rangle,
\end{eqnarray*}
where $v=e^{i\theta}\overline{Dh(z_0)}(z_0)+e^{-i\theta}h(z_0)$. This implies that $v$ is a normal vector of $\partial B(r)$ at $z_0$, thus there exists a real number $m$ such that $v=m z_0$. Since
\begin{eqnarray*}
m\|z_0\|^2&=&{\rm Re}\langle z_0, mz_0\rangle={\rm Re}\langle z_0, v\rangle\\
&=&{\rm Re}\ \langle z_0, e^{i\theta}\overline{Dh(z_0)}(z_0)\rangle+{\rm Re}\ \langle z_0,e^{-i\theta}h(z_0)\rangle  \\
&=&{\rm Re}\ [e^{-i\theta}\langle Dh(z_0)(z_0), z_0\rangle]+{\rm Re}\ \langle [e^{-i\theta}\langle h(z_0), z_0\rangle]\\
&=&  \frac{\|z_0\|^2}{2\rho}\ (\cos\theta +1 +t)+\frac{\|z_0\|^2}{2\rho}\ (\cos\theta+1),
\end{eqnarray*}
we obtain that $m=\frac{2\cos\theta + 2 + t}{2\rho}$, and this completes the proof of Lemma 2.

{\bf Proof of Theorem 4.}\ Let $h(z)=[Df(z)]^{-1}f(z)$ with $f(z)\in K(B, \alpha)$, and $g(z)=\frac{\|z\|^2}{\langle h(z), z\rangle}$, then $g(0)=1>\beta=\beta(\alpha)$.

Suppose $f\not\in S^{\ast}(B, \beta)$, then by the continuity of $g(z)$, there exists $z_0\in B\setminus\{0\}$ such that ${\rm Re}\ g(z_0)=\beta$ and ${\rm Re}\ g(z)\geq\beta$
for all $\|z\|<\|z_0\|$. Thus it follows from Lemma 2 that there exist real numbers $\theta$, $t\geq\frac{1-\beta}{\beta}$ and $m=\frac{2\cos\theta + 2 + t}{2\beta}$ such that (i)-(iii) of Lemma 2 hold.

Let $b=e^{-i\theta} h(z_0)-\frac{z_0}{2\beta} (1+\cos \theta)$, then it follows from Lemma 2(i) that
$$
{\rm Re}\ \langle b, z_0\rangle={\rm Re}\ [e^{-i\theta}\langle h(z_0), z_0\rangle ]-\frac{\|z_0\|^2}{2\beta} (1+\cos \theta)=0,
$$
so that
\begin{equation}
{\rm Re}\langle [Df(z_0)]^{-1}D^2f(z_0)(b,b), z_0\rangle < (1-\alpha)\ \|b\|^2.
\label{a4}
\end{equation}

Let $b_1=i z_0$, then ${\rm Re}\ \langle b_1, z_0\rangle={\rm Re}[i\|z_0\|^2]=0$, so we conclude from the fact $f(z)\in K(B, \alpha)$ that
\begin{equation}
{\rm Re}\langle [Df(z_0)]^{-1}D^2f(z_0)(iz_0,i z_0), z_0\rangle < (1-\alpha)\ \|z_0\|^2.
\label{a5}
\end{equation}

On the other hand, by Lemma 2, we have
\begin{equation}
\|b\|^2=\|h(z_0)\|^2-\frac{\|z_0\|^2}{4\beta^2} (1+\cos \theta)^2,
\label{a6}
\end{equation}
and $\|h(z_0)\|\geq\frac{\|z_0\|}{2\beta}|1+e^{i\theta}|=\frac{\|z_0\|}{2\beta}\sqrt{2+2\cos\theta}$, and
\begin{eqnarray*}
m\ e^{-i\theta}\langle h(z_0), z_0\rangle &=&\langle e^{-i\theta}h(z_0), mz_0\rangle\\
&=&\langle e^{-i\theta}h(z_0),e^{i\theta}\overline{Dh(z_0)}(z_0)\rangle+\langle e^{-i\theta}h(z_0), e^{-i\theta}h(z_0)\rangle \\
&=&e^{-i2\theta}\ \langle Dh(z_0)h(z_0), z_0\rangle+\|h(z_0)\|^2,
\end{eqnarray*}
so that
\begin{eqnarray*}
{\rm Re}\ [e^{-i2\theta}\ \langle Dh(z_0)h(z_0), z_0\rangle] &=&m{\rm Re}\ [e^{-i\theta}\langle h(z_0), z_0\rangle]-\|h(z_0)\|^2\\
&=&\frac{m\|z_0\|^2}{2\beta} (1+\cos \theta)-\|h(z_0)\|^2.
\end{eqnarray*}

By computing the Frechet derivatives for both sides of equation $Df(z)h(z)=f(z)$, we obtain
$$
[Df(z)]^{-1}D^2f(z)(h(z),h(z))+ Dh(z)h(z)=h(z),
$$
and
$$
[Df(z)]^{-1}D^2f(z)(h(z),z)+ Dh(z)(z)=z,
$$
thus we can obtain from the above equalities that
\begin{eqnarray*}
&&{\rm Re}\langle [Df(z_0)]^{-1}D^2f(z_0)(b,b), z_0\rangle ={\rm Re}[e^{-i2\theta}\langle [Df(z_0)]^{-1}D^2f(z_0)(h(z_0),h(z_0)), z_0\rangle \\
&&\hspace{3cm}-\frac{1+\cos\theta}{\beta}{\rm Re}[e^{-i\theta}\langle [Df(z_0)]^{-1}D^2f(z_0)(h(z_0),z_0), z_0\rangle ]\\
&&\hspace{3cm}+ \frac{(1+\cos\theta )^2}{4\beta^2}{\rm Re}[e^{-i\theta}\langle [Df(z_0)]^{-1}D^2f(z_0)(z_0,z_0), z_0\rangle ] \\
&&\hspace{3cm}=-{\rm Re}[e^{-i2\theta}\ \langle Dh(z_0)h(z_0), z_0\rangle ]+ {\rm Re}[e^{-i2\theta}\ \langle h(z_0), z_0\rangle ]\\
&&\hspace{3cm}+\frac{1+\cos\theta}{\beta}\{{\rm Re}[e^{-i\theta}\langle Dh(z_0)(z_0), z_0\rangle ]-{\rm Re}[e^{-i\theta}\langle z_0, z_0\rangle ]\}\\
&&\hspace{3cm}-\frac{(1+\cos\theta )^2}{4\beta^2}{\rm Re}[e^{-i\theta}\langle [Df(z_0)]^{-1}D^2f(z_0)(i z_0, i z_0), z_0\rangle ] \\
&&\hspace{3cm}\geq -\frac{m\|z_0\|^2}{2\beta} (1+\cos \theta)+\|h(z_0)\|^2-\frac{\|z_0\|^2}{2\beta} (1+\cos \theta)\\
&&\hspace{3cm}+\frac{\|z_0\|^2}{2\beta^2} (1+\cos \theta)(\cos\theta +1+t) +\frac{(1+\cos\theta)^2}{4\beta^2}(\alpha-1)\|z_0\|^2\\
&&\hspace{3cm}\geq(1-\alpha)\ \|b\|^2+\alpha\|h(z_0)\|^2-\alpha\frac{\|z_0\|^2}{4\beta^2}(2+2\cos\theta)\\
&&\hspace{3cm}+\frac{1+\cos\theta}{4\beta^2}\|z_0\|^2(2\alpha-2\beta+\frac{1-\beta}{\beta})\\
&&\hspace{3cm}\geq(1-\alpha)\ \|b\|^2,
\end{eqnarray*}
which contradicts (\ref{a4}). Hence $f\in S^{\ast}(B, \beta)$, and the proof is complete.

{\bf Remark 4.} Setting $X=\mathbb{C}^n$ in Theorem 4, we obtain the related result in \cite{k2, ll, z2}

By applying the growth theorem\cite{hgp1} of starlike mappings of order $\rho$ and Theorem 4, we have the following corollary.

{\bf Corollary 3.} Let $0\leq\alpha<1, \beta=\beta(\alpha)=\frac{2\alpha-1+\sqrt{(2\alpha-1)^2+8}}{4}$. If $f(z)\in K(B, \alpha)$, then for $z\in B$, we have
$$
\frac{\|z\|}{(1+\|z\|)^{2(1-\beta)}}\leq\|f(z)\|\leq\frac{\|z\|}{(1-\|z\|)^{2(1-\beta)}},
$$
and $f(B)\supset \frac{1}{2^{2(1-\beta)}}B$.

{\bf Remark 5.} Setting $\alpha=0$ in Corollary 3, we obtain the growth theorem of convex mappings\cite{gs,gk}.

Finally, we introduce a linear operator\cite{zl2} in purpose to construct some
concrete examples of biholomorphic convex mappings of order $\alpha$ on
$B$ in a Hilbert space $X$.

Let
$$
H(U)=\{f:U\rightarrow \mathbb{C} \mbox{ are analytic in } U
\mbox{ with } f(0)=0, f'(0)=1\},
$$
then
$$
f\in K(\alpha)\Longleftrightarrow f\in H(U )\mbox{ and }
\mbox{Re}\bigg\{1+\frac{zf''(z)}{f'(z)}\bigg\}>\alpha \quad\mbox{for
all} \quad z\in U.
$$
Let
$$
SK(B,\alpha)=\{f\in N(B): \|Df(z)^{-1}D^2f(z)(\cdot, \cdot)\|\leq 1-\alpha
\mbox{ for all } z\in B\}.
$$
From Theorem 1, we have  $SK(B,\alpha)\subset K(B,\alpha), SK(U,\alpha)\subset
K(\alpha)$ and
$$
SK(U, \alpha)=\bigg\{f\in H(U):
\bigg|\frac{f''(z)}{f'(z)}\bigg|\leq1-\alpha\mbox{ for all }
z\in U\bigg\}.
$$

Let $m$ be a positive integer and $\dim X\geq m\geq2$. Then there
exist $ u_1, u_2, \cdots, u_m\in X$ with $\|u_j\|=1 (j=1, 2,
\cdots, m)$ such that $\langle u_j, u_k\rangle=0 (j\neq k)$. For
$g_1, g_2, \cdots, g_m\in H(U)$, we define the operator
$\Phi$ as
\begin{eqnarray}
\Phi_{u_1, u_2, \cdots, u_m}(g_1, g_2, \cdots,
g_m)(z)=z-\he{j=1}{m}\langle z, u_j\rangle
u_j+\he{j=1}{m}g_j(\langle z, u_j\rangle)u_j \label{a21}
\end{eqnarray}
for $z\in B$.

{\bf Theorem 5.} Suppose that $0\leq\alpha<1, \Phi_{u_1, u_2, \cdots, u_m}(g_1,
g_2, \cdots, g_m)$ is defined by (\ref{a21}), where $g_1, g_2,
\cdots, g_m\in H(U)$ are locally univalent functions  on
$\Delta$.

(1) If $\Phi_{u_1,  \cdots, u_m}(g_1, \cdots, g_m)\in K(B,\alpha)$, then $g_1, g_2, \cdots, g_m\in K(\alpha)$.

(2) If $ h(\xi)=\left\{
\begin{array}{lll}
\frac{1-(1-\xi)^{2\alpha-1}}{2\alpha -1}, & \alpha\neq \frac{1}{2}\\
-\ln (1-\xi), & \alpha=\frac{1}{2}
\end{array}\right.$, then $h\in K(\alpha)$,
but $\Phi_{u_1, \cdots, u_m}(h, \cdots, h)\notin K(B,\alpha)$.

(3) $\Phi_{u_1, u_2, \cdots, u_m}(g_1, g_2, \cdots, g_m)\in SK(B,\alpha)$
if and only if $g_1, g_2, \cdots, g_m\in SK(U,\alpha)$.

{\bf Proof.} Let $f(z)=\Phi_{u_1, u_2, \cdots, u_m}(g_1, g_2,
\cdots, g_m)(z)$, where $g_1, g_2, \cdots, g_m\in H(U)$ are
locally univalent functions on $U$. By some straightforward
computations, we obtain
\begin{eqnarray*}
&&Df(z)=I-\he{j=1}{m}\langle\cdot, u_j\rangle u_j
+\he{j=1}{m}g'_j(\langle z, u_j\rangle)\langle\cdot, u_j\rangle u_j,\\
&&Df(z)^{-1}=I-\he{j=1}{m}\bigg(1-\frac{1}{g'_j(\langle z,
u_j\rangle)}\bigg)\langle\cdot, u_j\rangle u_j,\\
&&D^2f(z)(x, x)=\he{j=1}{m}g''_j(\langle z, u_j\rangle)[\langle x,
u_j\rangle]^2u_j
\end{eqnarray*}
for $z\in B$ and $x\in X$. Hence we have
\begin{eqnarray}
Df(z)^{-1}D^2f(z)(x, x)=\he{j=1}{m}\frac{g''_j(\langle z,
u_j\rangle)}{g'_j(\langle z, u_j\rangle)}[\langle x,
u_j\rangle]^2u_j. \label{a22}
\end{eqnarray}

(1) Assume that $f=\Phi_{u_1, u_2, \cdots, u_m}(g_1, g_2, \cdots,
g_m)\in K(B)$, for every $\xi\in U\setminus\{0\}$ and $k$
fixed, we let $z=\xi u_k$ and $x=i\xi u_k$, then $\mbox{Re}\langle
x, z\rangle=\mbox{Re}\{i|\xi|^2\}=0$ and $z\in B\setminus\{0\}$.
Note that $\langle u_j, u_k\rangle=0 (j\neq k)$, from (\ref{a22}),
we obtain
\begin{eqnarray*}
&&\|x\|^2-\mbox{Re}\langle Df(z)^{-1}D^2f(z)(x, x),
z\rangle\\
&=&|\xi|^2+|\xi|^2\mbox{Re}\bigg\{\frac{\xi g''_k(\xi)}{g'_k(\xi)}\bigg\}=|\xi|^2\mbox{Re}\bigg\{1+\frac{\xi
g''_k(\xi)}{g'_k(\xi)}\bigg\}>\alpha\|x\|^2=\alpha |\xi|^2
\end{eqnarray*}
for $\xi\in U\setminus\{0\}$. This implies $g_k\in K(\alpha)$ for
$k=1, 2, \cdots, m$.

(2) A simple computation yields
$$
\mbox{Re}\bigg\{1+\frac{\xi
h''(\xi)}{h'(\xi)}\bigg\}=\alpha+(1-\alpha)\mbox{Re}\bigg\{\frac{1+\xi}{1-\xi}\bigg\}>\alpha
$$
for all $\xi\in U$. It follows $h\in K(\alpha)$.

Let $\displaystyle\sqrt{\frac{1-\alpha}{2}}<a<1,
\displaystyle\frac{\sqrt{1-\alpha}}{\sqrt{2}a}<r<1$,
$$
z=ru_1+a\sqrt{1-r^2}u_2 \quad\mbox{and}\quad
x=a\sqrt{1-r^2}u_1-ru_2.
$$
Then we have
$\|x\|^2=a^2(1-r^2)+r^2>0, \mbox{Re}\langle x, z\rangle=0$ and
$$
0<\|z\|^2=r^2+a^2(1-r^2)<1.
$$
Notice that $\langle x, u_j\rangle=0(j\geq3)$, from (\ref{a22}) , we
obtain
\begin{eqnarray*}
\|x\|^2&-&\mbox{Re}\langle DF(z)^{-1}D^2F(z)(x, x),
z\rangle\\
&=&\|x\|^2-\mbox{Re}\bigg\{\frac{2}{1-\langle z,
u_1\rangle}[\langle x, u_1\rangle]^2\langle u_1,
z\rangle+\frac{2}{1-\langle z, u_2\rangle}[\langle x,
u_2\rangle]^2\langle u_2,
z\rangle\bigg\}\\
&=&a^2(1-r^2)+r^2-\bigg\{\frac{2r}{1-r}a^2(1-r^2)+\frac{2a\sqrt{1-r^2}}{1-a\sqrt{1-r^2}}r^2\bigg\}\\
&<&1-2r^2a^2<\alpha,
\end{eqnarray*}
where $F=\Phi_{u_1, u_2, \cdots, m}(h, h, \cdots, h)$. By Corollary 1, we have $F\not\in K(B, \alpha)$.

(3) Assume that $g_1, g_2, \cdots, g_m\in SK(U,\alpha), f=\Phi_{u_1,
u_2, \cdots, u_m}(g_1, g_2, \cdots, g_m)$, from (\ref{a22}), we
obtain
\begin{eqnarray}
\Big\|Df(z)^{-1}D^2f(z)(x,
x)\Big\|&=&\bigg\|\he{j=1}{m}\frac{g''_j(\langle z,
u_j\rangle)}{g'_j(\langle z, u_j\rangle)}[\langle x,
u_j\rangle]^2u_j\bigg\|\nonumber\\
&\leq&\he{j=1}{m}\bigg|\frac{g''_j(\langle z,
u_j\rangle)}{g'_j(\langle z, u_j\rangle)}\bigg||\langle x,
u_j\rangle|^2\nonumber\\
&\leq&(1-\alpha)\he{j=1}{m}|\langle x, u_j\rangle|^2\label{a24}
\end{eqnarray}
for $z\in B$ and $x\in X$.

Fix $x\in X$, let $x_0=\he{j=1}{m}\langle x, u_j\rangle u_j$, a
simple computation yields
$$
\langle x-x_0, u_j\rangle=\langle x, u_j\rangle-\he{k=1}{m}\langle
x, u_k\rangle\langle u_k, u_j\rangle=\langle x, u_j\rangle-\langle
x, u_j\rangle=0,
$$
for $j=1, 2, \cdots, m$. This leads to $\langle x-x_0,
x_0\rangle=0$. Hence we conclude that
\begin{eqnarray}
\|x\|^2&=&\|(x-x_0)+x_0\|^2=\|x-x_0\|^2+\|x_0\|^2\nonumber\\
&=&\|x-x_0\|^2+\he{j=1}{m}|\langle
 x, u_j\rangle|^2\nonumber\\
 &\geq&\he{j=1}{m}|\langle
 x, u_j\rangle|^2. \label{a25}
\end{eqnarray}
From (\ref{a24}) and (\ref{a25}), we obtain
$$
\Big\|Df(z)^{-1}D^2f(z)(x, x)\Big\|\leq (1-\alpha)\he{j=1}{m}|\langle x,
u_j\rangle|^2\leq (1-\alpha)\|x\|^2\leq 1-\alpha
$$
for all $z\in B$ and $x\in X$ with $\|x\|=1$. Since $X$ is a
Hilbert space, by the result in \cite{tl}(see P.342), we have
\begin{eqnarray*}
\|Df(z)^{-1}D^2f(z)(\cdot, \cdot)\|
&=&\sup_{\|x\|=1, \ \|y\|=1}\|Df(z)^{-1}D^2f(z)(x, y)\|\\
&=&\sup_{\|x\|=1}\|Df(z)^{-1}D^2f(z)(x, x)\| \leq 1-\alpha.
\end{eqnarray*}
It follows that $f\in SK(B,\alpha)$.

Conversely, suppose that $f=\Phi_{u_1, u_2, \cdots, u_m}(g_1, g_2,
\cdots, g_m)\in SK(B,\alpha)$.  For every $\xi\in U$ and $k$ fixed
($1\leq k\leq m$), we let $z=\xi u_k$ and $x=u_k$, then we have
$z\in B, \langle z, u_k\rangle=\xi$ and $\|x\|=1$. Note that
$\langle u_j, u_k\rangle=0 (j\neq k)$ and $\|u_k\|=1$, from
(\ref{a22}), we obtain
$$
\bigg|\frac{g''_k(\xi)}{g'_k(\xi)}\bigg|=\|Df(z)^{-1}D^2f(z)(x,
x)\|\leq\|Df(z)^{-1}D^2f(z)(\cdot, \cdot)\|\|x\|^2\leq 1-\alpha
$$
for $\xi\in U$. That is, $g_k\in SK(U, \alpha)$ for $k=1, 2,
\cdots, m$, and the proof is complete.

{\bf Remark 6.} Let $X=\mathbb{C}^n, \alpha=0$. If we choose $u_1=(1, 0, \cdots, 0),
u_2=(0, 1, \cdots, 0), \cdots$, $u_n=(0, 0, \cdots, 1)\in \mathbb{C}^n$,
then we have $z=\he{j=1}{n}\langle z, u_j\rangle u_j$ for $z\in
\mathbb{C}^n$. From Theorem 5, we obtain a result, which is
Theorem 3 and Theorem 4 in \cite{z} for case  $p=2$. Part (2) was
obtained by Roper and Suffridge \cite{rs1}, \cite{rs} using a
different method.

{\bf Example 3.} Let $0\leq\alpha<1, \dim X\geq m\geq2$ and $\lambda_j\in \mathbb{C}$ with
$\lambda_j\neq 0(j=1, 2, \cdots, m)$, then
$$
f(z)=z-\he{j=1}{m}\langle z, u_j\rangle
u_j+\he{j=1}{m}\frac{e^{\lambda_j \langle z,
u_j\rangle}-1}{\lambda_j}u_j\in K(B, \alpha)\Longleftrightarrow
|\lambda_j|\leq1-\alpha (j=1, 2, \cdots, m),
$$
where $ u_j\in X$ with $\|u_j\|=1 $ such that $\langle u_j,
u_k\rangle=0 (j, k=1, 2, \cdots, m, j\neq k)$.

{\bf Proof.} If $|\lambda_j|\leq1(j=1, 2, \cdots, m)$,  setting
$g_j(\xi)=\displaystyle\frac{e^{\lambda_j\xi}-1}{\lambda_j}$, then
we have that $g_j$ is analytic on $U$ with $g_j(0)=0,
g'_j(0)=1$ such that $
\displaystyle\frac{g''_j(\xi)}{g'_j(\xi)}=\lambda_j $ for
$\xi\in U(j=1, 2, \cdots, m)$. Hence $g_j\in SK(U)$. From
Theorem 5, we obtain $f\in SK(B, \alpha)$.

Conversely, we shall prove that $|\lambda_j|\leq 1-\alpha$ for all $j=1,
2, \cdots, m$ when $f$ is a biholomorphic convex mapping of order $\alpha$ on $B$.

If not, then there exists $k$ such that $|\lambda_k|>1-\alpha$. Let
$\displaystyle\frac{1-\alpha}{|\lambda_k|}<r<1, \theta=\arg\lambda_k,
z_0=-re^{-i\theta}u_k$ and $x=ie^{-i\theta}u_k$, then $\|x\|=1,
\mbox{Re}\langle x, z_0\rangle=\mbox{Re}\{-ir\}=0$. Using the fact
that $\langle u_j, u_k\rangle=0(j\neq k)$, from (\ref{a22}), we
obtain
\begin{eqnarray*}
Df(z_0)^{-1}D^2f(z_0)(x, x)&=&\he{j=1}{m}\frac{g''_j(\langle z_0,
u_j\rangle)}{g'_j(\langle z_0, u_j\rangle)}[\langle x,
u_j\rangle]^2u_j\\
&=&\lambda_k[\langle x,
u_k\rangle]^2u_k=-|\lambda_k|e^{-i\theta}u_k.
\end{eqnarray*}
Hence we have
$$
\|x\|^2-\mbox{Re}\langle Df(z_0)^{-1}D^2f(z_0)(x, x),
z_0\rangle=1-r|\lambda_k|<\alpha,
$$
which contradicts (\ref{a02}). This completes the proof.

\end{document}